\documentclass[a4paper,10pt]{article}  
\usepackage{amssymb}

\usepackage{graphicx}
\usepackage{psfrag}
\usepackage{graphicx}
\usepackage{amsmath}
\usepackage{breqn}
\usepackage{ulem}
\newcommand{\disp}{\displaystyle}
\usepackage{latexsym}
\usepackage{xcolor}

\newtheorem{theorem}{Theorem}[section]
\newtheorem{corollary}[theorem]{Corollary}
\newtheorem{lemma}[theorem]{Lemma}
\newtheorem{proposition}[theorem]{Proposition}

\newtheorem{definition}[theorem]{Definition}

\numberwithin{equation}{section}

\begin{document}

\title{A Maximum Principle for Optimal Control Problems involving Sweeping Processes with a Nonsmooth  Set}

\author{ M. d. R.   de Pinho,  M. Margarida A. Ferreira \thanks{MdR de Pinho  and MMA Ferreira  are at Faculdade de Engenharia da Universidade do Porto, DEEC, SYSTEC. 
              Portugal,            
             mrpinho, mmf@fe.up.pt}  and  \hspace{0.1cm} Georgi Smirnov  \thanks{G. Smirnov  is at
            Universidade do Minho, Dep. Matem\'{a}tica,
            Physics Center of Minho and Porto Universities (CF-UM-UP), Campus de Gualtar, Braga, Portugal,         
             smirnov@math.uminho.pt}
             }

\maketitle

\begin{abstract}
We generalize a  Maximum Principle for  optimal control problems involving sweeping systems  previously derived in   \cite{nosso_2022} to cover the case where the moving set may be nonsmooth. Noteworthy, we consider problems with constrained end point.  A remarkable feature of our work is that we rely upon an ingenious  smooth approximating family of standard differential equations in the vein of that used in \cite{nosso_2019}.
\end{abstract}


%All acknowledgements should be placed in the back of the paper after Conclusions..

\section{Introduction}
%Sweeping  processes are evolution differential inclusions involving the normal cone to a set.  Introduced   in the seminal paper  \cite{Mo74}  by J.J. Moreau in  the context of  plasticity and friction theory,  these systems have proved of interest  to tackle  problems  in   mechanics, engineering, economics  and  crowd motion problems; see, for example,   \cite{Addy}, \cite{MoCa17}, \cite{KuMa00}, \cite{Maury},  \cite{Thib2016} and \cite{nosso_2019}. Recently, there has been a surge of interest in 
%optimal control problem involving controlled sweeping systems of the form
In recent years,  there has been a surge of interest in optimal control problems involving the controlled sweeping process of the form
 \begin{equation}\label{SP}
\dot x(t) \in f(t,x(t),u(t))- N_{C(t)}(x(t)), ~u(t)\in U, ~~x(0) \in C_0.
\end{equation} 
In this respect,  we refer to, for example, 
 \cite{ArCo17},  \cite{BrKr},  \cite{MoCa17}, \cite{CoPa16},    \cite{CoHeHoMo}, \cite{KuMa00}, \cite{zeidan2020}, \cite{nosso_2019} (see also accompanying correction \cite{correction_2019}), \cite{CCBN_2021},  \cite{Palladino2022} and \cite{nosso_2022}.  
 Sweeping processes first appeared in the seminal paper  \cite{Mo74}  by J.J. Moreau as a mathematical framework for problems in   plasticity and friction theory. They  have proved of interest  to tackle  problems  in   mechanics, engineering, economics  and  crowd motion problems; to name but a few, see \cite{Addy}, \cite{MoCa17}, \cite{KuMa00}, \cite{Maury} and  \cite{Thib2016}.
In the last decades, systems  in the form  (\ref{SP})  have  caught the attention and interest of the optimal control community. Such interest resides not only  in the range of applications but  also in the   remarkable challenge they rise concerning  the derivation of necessary conditions. This is due to   the  presence of the normal cone $N_{C(t)}(x(t))$  in  the dynamics. Indeed, the presence of the normal cone renders the discontinuity of the right hand of the differential inclusion in   (\ref{SP})  destroying a regularity property central  to many  known  optimal control results.

Lately, there has been several successful attempts to derive necessary conditions for optimal control problems involving  (\ref{SP}). Assuming that the set $C$ is time independent,  necessary conditions  for  optimal control  problems with free end point   have been derived under different assumptions and using different techniques. 
In  \cite{nosso_2019}, the set $C$ has the form $C=\{ x :~\psi(x)\leq 0\}$ and an approximating sequence  of optimal control problems, where   (\ref{SP})  is  approximated  by  the    differential equation  
\begin{equation}\label{Sgamma}
\dot x_{\gamma_k}(t)=  f(t,x_{\gamma_k}(t),u(t))-\gamma_k e^{\gamma_k \psi(x_{\gamma_k}(t))}\nabla \psi(x_{\gamma_k}(t)), \end{equation}
for some positive sequence $\gamma_k\to +\infty$, is used. Similar techniques are also applied  to somehow more general problems in \cite{zeidan2020}. 
A useful feature of those approximations   is explored in \cite{nosso_Jota} to define numerial schemes to solve such problems.

More recently,  an adaptation of the family of approximating  systems \eqref{Sgamma}   is used in \cite{nosso_2022} to generalize the results in \cite{nosso_2019} to  cover problems with additional end point constraints and   with a moving set of the form  $C(t)=\{ x :~\psi(t,x)\leq 0\}$.   

In this paper we generalize  the Maximum Principle proved in \cite{nosso_2022}   to cover problems with possibly nonsmooth sets. Our problem of interest is
$$
(P)
\left\{
\begin{array}{l}
\mbox{Minimize } \; 
\phi(x(T))\\[2mm]
\mbox{over processes $(x,u)$ such that }\\[2mm]
\hspace{8mm} \dot x(t) \in f(t,x(t),u(t))- N_{C(t)}(x(t)), \hspace{0.2cm}\mbox{a.e.}\ \ t\in[0,T],\\[2mm]
\hspace{8mm} u(t)\in U, \ \ \;\, \mbox{a.e.}\ \ t\in[0,T],\\[2mm]
\hspace{8mm} (x(0),x(T)) \in C_0\times C_T { ~\subset C(0)\times C(T)},
\hspace{2mm}\end{array}
\right.
$$
where  $T>0$ is fixed, $\phi: R^n\to  R$, $f:[0,T]\times R^n\times  R^m\to  R^n$, $U \subset  R^m$ %$C_0\subset C(0),~C_T \subset C(T)$,
 and
\begin{equation}
\label{set:C}
C(t):=\left\lbrace x\in  R^n: ~\psi^i(t,x)\leq 0,\; i=1,\ldots, I\right\rbrace
\end{equation}
for some  functions $\psi^i:[0,T]\times  R^n\to  R$, $ i=1,\ldots, I$.

The  case where   $I=1$ in \eqref{set:C}  and $\psi^1$ is $C^2$ is covered in \cite{nosso_2022}. Here, we  assume $I>1$ and that the functions $\psi^i$ are also $C^2$.
Although going from $I=1$ in \eqref{set:C} to $I>1$ may be seen as a small generalization, it demands a significant revision of the technical approach and, plus,  the introduction of a constraint qualification. This is because   the set  \eqref{set:C}  may be nonsmooth.  We focus on sets  \eqref{set:C},  satisfying a certain constraint qualification,  introduced in assumption (A1) in section \ref{Prel} below. This is, indeed, a  restriction on the nonsmoothness of \eqref{set:C}.
A similar problem with nonsmooth moving set is considered in \cite{Palladino2022}. Our results cannot be obtained from the results of \cite{Palladino2022} and do not generalize them.

This paper is organized in the following way. In section \ref{Prel}, we introduce the main notation and we state and discuss the assumptions under which we work. In this same section, we also introduce the family of approximating systems to $\dot x(t) \in f(t,x(t),u(t))- N_{C(t)}(x(t))$ and establish a crucial convergence result, Theorem  \ref{T1}. In section \ref{ApFam}, we dwell on the approximating family of optimal control problems to $(P)$ and we state the associated necessary conditions. The Maximum Principle for $(P)$ is then deduced and stated in Theorem \ref{main}, covering additionally, problems in the form of $(P)$ where the end point constraint $x(T) \in  C_T$ is absent. 
Before finishing, we present   an illustrative example of our main result, Theorem \ref{main}. 
%A section with conclusions ends this paper.

\section{Preliminaries}
\label{Prel}
In this section, we introduce a summary of the notation and  state  the assumptions  on the data of $(P)$ enforced throughout.  Furthermore, we extract information from the assumptions  establishing relations  crucial for  the forthcoming analysis.
\medskip

\noindent {\bf Notation}

For a set $S\subset  R^n$, $\partial S$, $  \text{cl}\,S$ and $  \text {int}\, S$ denote the \textit{boundary}, \textit{closure}   and   \textit{interior} of $S$. 

If  $g: R^p\to  R^q$, $\nabla g$ represents the derivative and $\nabla^2g$ the second derivative.   
If $g: R\times  R^p\to  R^q$,  then $\nabla_x g$ represents the derivative w.r.t. $x\in  R^p$ and $\nabla^2_xg$ the second derivative, while $\partial_t g(t,x)$ represents the derivative w.r.t. $t\in  R$.   

The Euclidean norm or the induced matrix
norm on $ R^{p\times q}$ is denoted by $  |\cdot |$.  We denote by $B_n$ the closed unit ball in $ R^n$  centered at the origin. The inner product  of $x$ and $y$ is denoted by $\langle x, y\rangle $. For some $A\subset  R^n$, $d(x,A)$ denotes the distance between $x$ and $A$.  We denote the support function of $A$ at $z$ by $S(z,A)=\sup\{\langle z,a\rangle\mid a\in A\}$ 

The space  $L^{\infty}([a,b]; R^p)$ (or simply $L^{\infty}$ when the domains are clearly understood) is the Lebesgue space of essentially 
bounded functions  $h:[a,b]\to R^p$. We say that $h\in BV([a,b]; R^p)$ if $h$ is a function of bounded variation. The space of continuous functions is denoted by $C([a,b]; R^p)$.
%The norm of $L^\infty([a,b]; R^p)$ or $C([a,b]; R^p)$ is denoted by $\|\cdot\|_\infty $. 
%Let $C^*([a,b]; R)$ 
%be the dual space of the continuous functions defined from $[a,b]$ to $ R$, denoted by 
%$C([a,b]; R)$, with supremum norm. The norm of $C^*([a,b]; R)$ is denoted by $\|\mu\|_{TV}$. We denote $C^\oplus([a,b]; R)$ the set of elements in $C^*([a,b]; R)$  which take nonnegative values on nonnegative valued functions in $C([a,b]; R)$. For $\mu \in C^\oplus([a,b]; R)$,
%$\|\mu\|_{TV}=\int_{[a,b]} \mu(dt)$.

Standard concepts from nonsmooth analysis will also be used. Those can be found in  \cite{Cla83}, \cite{Mor05} or \cite{Vinter}, to name but a few. 
The \textit{Mordukhovich} normal cone to a set $S$  at $s\in S$  is denoted by 
$N_{S}(s)$   and 
 $\partial f(s)$ is the
\textit{Mordukhovich} subdifferential of $f$ at $s$ (also known as \textit{limiting subdifferential}). 

For any set $A\subset  R^n$, $  \text{cone}\, A$ is the cone generated by the set $A$.
\bigskip

We now turn to problem $(P)$.  We first state    the definition of admissible processes for $(P)$  and then we describe the assumptions under which we will derive our main results.

\begin{definition}
A pair  $(x,u)$ is called an admissible process for $(P)$ when  $x$ is an absolutely continuous function and  $u$ is a measurable function  satisfying the constraints of $(P)$.
\end{definition}

\noindent {\bf Assumptions on the data of $\mathbf{(P)}$}
\begin{itemize}
\item[A1:]  The function $\psi^i$, $ i=1,\ldots, I$, are  $C^2$.  The graph of $C(\cdot)$ is compact and   it  is  contained in the interior of a ball $rB_{n+1}$, for some $r>0$. 
There exist constants $\beta>0$,  $\eta >0$ and $\rho \in ]0,1[$   such that
\begin{equation}
\psi^i(t,x) \in [ -\beta,\beta] \Longrightarrow |\nabla_x \psi^i (t,x) | > \eta\; \;   {\rm for all }\; (t,x)\in [0,T]\times   R^n,  \label{psi:bder}
\end{equation}
%\begin{equation}
%\psi^i(t,x) \leq -2\beta \Longrightarrow \nabla_x \psi^i (t,x)=0\; \;   {\rm for all }\; (t,x)\in [0,T]\times   R^n,  \label{psi:bder1}
%\end{equation}
and, for  $I(t,x)=\{  i=1,\ldots, I\mid \psi^i(t,x)\in ]-2\beta,\beta]\}$,  
\begin{equation}
\label{cn1}
 \langle\nabla_x\psi^i(t,x),\nabla_x\psi^j(t,x)\rangle\geq 0,\;\; i,j\in I(t,x).
\end{equation}
Moreover,  if $i\in I(t,x)$, then
\begin{equation}
\label{cn2}
\sum_{j\in I(t,x)\setminus\{i\}} \big| \langle\nabla_x\psi^i(t,x),\nabla_x\psi^j(t,x)\rangle\big| \leq \rho |\nabla_x\psi^i(t,x)|^2
\end{equation}
and 
\begin{equation}\label{gradient0}
\psi^i(t,x)\leq -2\beta~\Longrightarrow~\nabla \psi^i(t,x)=0~\text{ for } i=1,\ldots I.
\end{equation}

\item[A2:] The function $f$ is  continuous, $x\to f(t,x,u)$  is continuously differentiable  for all $(t,u)\in [0,T]\times  R^m$. The constant  $M>0$ is such that $|f(t,x,u)|\leq M$ and $|\nabla_x f(t,x,u)|\leq M$ for all $(t,x,u)\in rB_{n+1}\times U$.
\item[A3:]  For each $(t,x)$,  the set $f(t,x,U)$ is  convex. 

\item[A4:] The set $U$ is compact.
\item[A5:] The sets  $C_0$ and $C_T$ are  compact.
\item[A6:] There exists a constant $L_\phi$ such that $|\phi(x)-\phi(x')|\leq L_{\phi}|x-x'|$ for all $x, x' \in  R^n$.
\end{itemize}
 
\begin{figure}[h!]
\begin{center}
\includegraphics{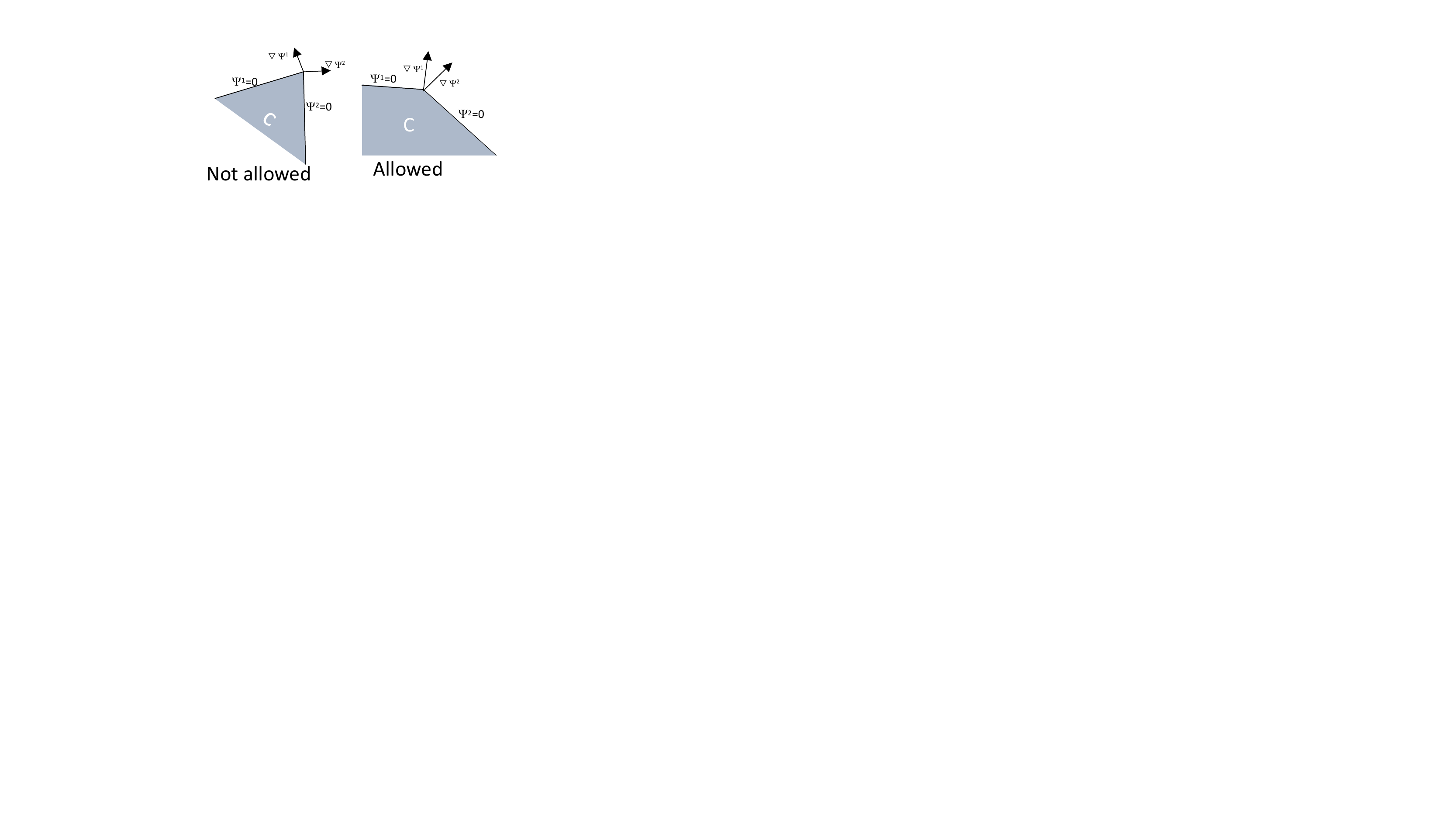}
\caption{Examples of two diferent sets $C$. On the left size, a set  that does not satisfies  \eqref{cn1}. On the right side, the set $C$ is nonsmooth and it fulfils  \eqref{cn1}.}
\label{fig1}
\end{center}
\end{figure}

Assumption (A1) concerns the functions $\psi^i$ defining the set $C$ and it  plays a crucial role in the analysis. All   $\psi^i$ are assumed to be smooth with gradients   bounded away from the origin when $\psi^i$ takes values in a neighorhood of zero.  Moreover,  the boundary of $C$ may be  nonsmooth at the  intersection points  of the level sets  $\left\{x: \psi^i(t,x)=0\right\}$. However, nonsmoothness at those  corner points is restricted to  (\ref{cn1})  which excludes the cases where   the angle between the  two gradients of the functions defining the boundary of $C$  is obtuse; see figure \ref{fig1}.

On the other hand,  \eqref{cn2} guarantees that the Gramian matrix of the gradients of the functions taking values near the boundary of $C(t)$ is diagonally dominant and, hence,  the gradients are linearly independent. 

In many situations, as in the example we present in the last section, we can guarantee the fulfillment of (A1), in particular \eqref{gradient0}, replacing the function $\psi^i$ by
\begin{equation}\label{tilde:psi}\tilde \psi^i(t,x)=h\circ \psi^i(t,x), 
\end{equation}
where 
$$h(z)=\left\{\begin{array}{lcl}
z & \text{ if} & z>-\beta, \\
h_s (z) &  \text{ if} &-2\beta\leq  z\leq -\beta,\\
-2\beta   & \text{ if} & z<-2\beta,
\end{array}\right.$$
Here, $h$  is an  $C^2$  function, with  $h_s$ an  increasing function defined on  $[-2\beta,-\beta]$. For example, $h$ may be a cubic polinomial  with positive  derivative on the interval  $]-2\beta,-\beta[$. For all $t\in [0,T]$, set  $$\tilde C(t): =\left\lbrace x\in R:~\tilde \psi^i(t,x )\leq 0,~i=1, \ldots,I\right\rbrace.$$
It is then a simple matter to see that 
$$C(t)=\tilde C(t)  \text{ for all } t \in [0,T].$$
and that  the functions $\tilde \psi^i(\cdot)$ satisfy the assumption  (A1). 
%It follows that all the previous analysis  holds when we replace $\psi^i$ by $\tilde \psi^i$, for $i=1, \ldots I$. The replacement of $\psi^i$ by $\tilde \psi^i$ is of interest since the following property holds for $\tilde \psi^i$ and does not necessarily holds for $\psi^i$:
%\begin{equation}
%\tilde \psi^i(t,x) \leq -2\beta \Longrightarrow \nabla_x \tilde  \psi^i (t,x)=0\; \;   {\rm for~ all }\; (t,x)\in [0,T]\times   R^n.  \label{psi:bder1}
%\end{equation}

The assumption that the graph of $C(\cdot)$ is  compact  and contained in the interior of a ball  is introduced to avoid technicalities in our forthcoming analysis.  In applied problems,  this may be easily side tracked  by considering the intersection of the graph of   $C(\cdot)$ with a tube around the optimal trajectory.

\bigskip

We now proceed introducing  an approximation family of controlled systems  to \eqref{SP}. 
Let $x(\cdot)$ be a solution to  
 the differential inclusion
$$ \dot x(t) \in f(t,x(t),U)- N_{C(t)}(x(t)).$$
Under our assumptions, measurable selection theorems assert the existence of measurable functions $u$ and $\xi^i$ such that 
$u(t) \in U$, $\xi^i(t)\geq 0$ a.e.  $t\in [0,T]$, $\xi^i(t)=0$ if $\psi^i(t,x(t))<0$, and 
$$
\dot x(t)= f(t,x(t),u(t))-\sum_{i=1}^I\xi^i(t)\nabla_x\psi^i(t,x(t))\;  {\rm a.e.}\; t\in [0,T].
$$

Considering the   trajectory $x$, some  observations are called for.
Let $\mu$ be such that
\begin{multline}\nonumber
\max\left\{(|\nabla_x\psi^i(t,x)||f(t,x,u)|+|\partial_t\psi^i(t,x)|)+1:\right. \\[1mm]
\left.  ~t\in [0,T],\; u\in U,\; x\in C(t)+B_n, \; i=1,\ldots,I\right\}\leq\mu.
\end{multline}
The properties of the graph of $C(\cdot)$ in (A1) guarantee the existence of such  maximum.

Consider now some  $t$    such that, for some $j\in \{1, \ldots I\}$, $\psi^j (t,x(t))=0$ and $\dot{x}(t)$ exists. Since the trajectory $x$ is always in $C$, we have (see \eqref{cn1})
\begin{equation}\nonumber
\begin{split}
0 & =\frac{d}{dt}\psi^j(t,x(t))=\langle \nabla_x\psi^j(t,x(t)),\dot{x}(t)\rangle +\partial_t\psi^j(t,x(t))
\\[2mm]
& =\langle \nabla_x\psi^j(t,x(t)),f(t,x(t),u(t))\rangle-\xi^j(t)| \nabla_x\psi^j(t,x(t))|^2 
\\[2mm]
&
-\sum_{i\in I(t,x(t))\setminus\{ j\}}\xi^i(t) \langle\nabla_x\psi^i(t,x(t)),\nabla_x\psi^j(t,x(t))\rangle
+\partial_t\psi^j(t,x(t))
\\[2mm]
&
\leq
\langle \nabla_x\psi^j(t,x(t)),f(t,x(t),u(t))\rangle-\xi^j(t)| \nabla_x\psi^j(t,x(t))|^2 
+\partial_t\psi^j(t,x(t)),
\end{split}
\end{equation}
and, hence (see (\ref{psi:bder})),
$$
\xi^j(t)\leq\frac{1}{| \nabla_x\psi^j(t,x(t))|^2}(\langle \nabla_x\psi^j(t,x(t)),f(t,x(t),u(t))\rangle+\partial_t\psi^j(t,x(t)))
\leq \frac{\mu}{\eta^2}.
$$

Define the function 
$$\mu(\gamma)=\frac{1}{\gamma}\log\left( \frac{\mu}{\eta^2\gamma}\right), \quad \gamma>0,$$
%Let $\{\gamma_k\}$  be a sequence of positive scalars such that $\gamma_k\uparrow\infty$. Additionally, consider a  sequence  $\{\mu_k\}$ such that
%$$
%\mu_k=\frac{1}{\gamma_k}\log\frac{\mu}{\eta^2\gamma_k}>-\beta,
%$$
% for  $k$ is sufficiently large. 
% $\sigma>0$.%and define  the sets
% $$C_{(\sigma,\mu)}(t)=\left\{x: \psi(t,x(t))-\sigma\leq \mu\right\}.$$
consider a sequence $\{\sigma_k\}$ such that $\sigma_k\downarrow 0$ and choose another  sequence  $\{\gamma_k\}$   with $\gamma_k\uparrow +\infty $ and 
$$C(t)\subset   {\rm int}\: C^k(t)=   {\rm int}\: \left\{x: \psi^i(t,x)-\sigma_k\leq \mu_k,\; i=1,\ldots,I\right\}, $$
where 
$$\mu_k=\mu(\gamma_k).$$
% Define
%$$
%C^k(t)=\{ x\mid \psi(t,x)\leq \mu_k\}
%$$
Let   $x_k$  be a solution to the differential equation
\begin{equation}\label{Approx}
\dot{x}_k(t)=f(t,x_k(t),u_k(t))-\sum_{i=1}^I\gamma_k e^{\gamma_k(\psi^i(t,x_k(t))-\sigma_k)}\nabla_x\psi^i(t,x_k(t))
\end{equation}
for some $u_k(t) \in U$ a.e. $t\in [0,T]$.
Take any  $t\in [0,T]$ such that $\dot{x}_k(t)$ exists and $\psi^j (t,x_k(t))-\sigma_k=\mu_k$. Assume $k$   is such that    $j\in I(t,x_k(t))$. Then,  whenever $\gamma_k$ is sufficiently large, we have
\begin{equation}\nonumber 
\begin{split}
\frac{d}{dt}\psi^j(t,x_k(t)) & =\langle \nabla_x\psi^j(t,x_k(t)),f(t,x_k(t),u_k(t))\rangle\\
& \hspace{1cm}
-\gamma_ke^{\gamma_k(\psi^j(t,x_k(t))-\sigma_k)}|\nabla_x\psi^j(t,x_k(t))|^2
\\
&  \hspace{0.3cm}
-\sum_{i\in I(t,x_k(t))\setminus\{ j\}}\hspace{-2em} \gamma_ke^{\gamma_k(\psi^i(t,x_k(t))-\sigma_k)}\langle\nabla_x\psi^i(t,x_k(t)),\nabla_x\psi^j(t,x_k(t))\rangle\\
& \hspace{0.3cm}
-\sum_{i\not\in I(t,x_k(t))}\hspace{-1em} \gamma_ke^{\gamma_k(\psi^i(t,x_k(t))-\sigma_k)}\langle\nabla_x\psi^i(t,x_k(t)),\nabla_x\psi^j(t,x_k(t))\rangle\\ &   \hspace{0.3cm} +
\partial_t\psi^j(t,x_k(t))\\
& 
\leq
\langle \nabla_x\psi^j(t,x_k(t)),f(t,x_k(t),u_k(t))\rangle\\
& \hspace{0.3cm} -\gamma_ke^{\gamma_k(\psi^j(t,x_k(t))-\sigma_k)}|\nabla_x\psi^j(t,x_k(t))|^2\\
& \hspace{0.3cm}
-\sum_{i\not\in I(t, x_k(t))} \hspace{-1em}\gamma_ke^{\gamma_k(\psi^i(t,x_k(t))-\sigma_k)}\langle\nabla_x\psi^i(t,x_k(t)),\nabla_x\psi^j(t,x_k(t))\rangle\\
& \hspace{0.3cm}+
\partial_t\psi^j(t,x_k(t))
\\ & 
\leq
\langle \nabla_x\psi^j(t,x_k(t)),f(t,x_k(t),u_k(t))\rangle\\
& \hspace{0.3cm} -\gamma_ke^{\gamma_k(\psi^j(t,x_k(t))-\sigma_k)}|\nabla_x\psi^j(t,x_k(t))|^2 \\
&  \hspace{0.3cm}
+\sum_{i\not\in I(t,x_k(t))}\hspace{-1em} \gamma_ke^{\gamma_k(-{\color{red}2}\beta-\sigma_k)}|\langle\nabla_x\psi^i(t,x_k(t)),\nabla_x\psi^j(t,x_k(t))\rangle|
\\
& \hspace{0.3cm}
+
\partial_t\psi^j(t,x_k(t))\\
& 
\leq \mu-\frac{1}{2} -\eta^2\gamma_k e^{\gamma_k\mu_k}\\
& =-\frac{1}{2}.
\end{split}
\end{equation}

Above,  we have used the definition of $\mu$ and the inequality 
$$\sum_{i\not\in I(t,x_k(t))}\hspace{-1em} \gamma_ke^{\gamma_k(-{\color{red}2}\beta-\sigma_k)}|\langle\nabla_x\psi^i(t,x_k(t)),\nabla_x\psi^j(t,x_k(t))\rangle|\leq \frac{1}{2},$$
which holds for $\gamma_k$ sufficiently large. 

Now, if  $x_k(0)\in C^k(0)$,  we  assure that  $x_k(t)\in C^k(t)$, for all $t\in [0,T]$, and
\begin{equation}
\gamma_ke^{\gamma_k(\psi^j(t,x_k(t))-\sigma_k)}\leq \gamma_ke^{\gamma_k\mu_k}= \frac{\mu}{\eta^2}.
\label{star}
\end{equation}
It follows that, for  $k$ sufficienttly large, we have
$$
|\dot{x}_k(t)|\leq ({\rm const}).
$$

 We are now a in position to state and prove our first result, Theorem \ref{T1} below. This is in the vein of  
 Theorem 4.1 in \cite{zeidan2020} (see also Lemma 1 in \cite{nosso_2019} when $\psi$ is independent of $t$ and convex) deviating from it in so far as the   approximating sequence of control systems \eqref{Approx} differs from the one introduced in  \cite{nosso_2019}\footnote{See also Theorem 2.2 in \cite{nosso_2022}}.  The proof of Theorem \ref{T1} relies on  (\ref{star}).
% ISTO ESTÀ COMO O ANTERIORThe choice of sequences $\mu_k$ and $\gamma_k$ above guarantees the  estimation (\ref{star}),   greatly simplifying  the proof of Theorem \ref{T1}. 
\begin{theorem}
\label{T1}
Let $\{(x_k,u_k)\}$, with $u_k(t)\in U$ a.e., be a sequence of solutions of Cauchy problems
\begin{equation}
\label{e1}
\left\{\begin{array}{rcl}
\dot{x}_k(t) & = & f(t,x_k(t),u_k(t))-\displaystyle\sum_{i=1}^I\gamma_k e^{\gamma_k(\psi^i(t,x_k(t))-\sigma_k)}\nabla_x\psi^i(t,x_k(t)),\\[2mm]
 x_k(0) & = & b_k\in C^k(0).
 \end{array}
 \right.
\end{equation}
If $b_k\rightarrow x_0$, then there exists a subsequence $\{x_k\}$ (we do not relabel) converging uniformly  to $x$, a unique solution
to the Cauchy problem
\begin{equation}
\label{e2}
\dot{x}(t)\in f(t,x(t),u(t))-N_{C(t)}(x(t)),\;\;\; x(0)=x_0,
\end{equation}
where $u$ is a  measurable function such that $u(t)\in U$ a.e. $t\in [0,T]$.

If, moreover, all the controls $u_k$ are equal, i.e., $u_k=u$, then the subsequence converges to a unique solution of (\ref{e2}), i.e., any solution of 
\begin{equation}
\label{e3}
\dot{x}(t)\in f(t,x(t),U)-N_{C(t)}(x(t)),\;\;\; x(0)=x_0\in C(0)
\end{equation}
can be approximated by solutions of (\ref{e1}).
\end{theorem}

\textit{ Proof}
Consider the sequence $\{x_k\}$, where $(x_k,u_k)$ solves (\ref{e1}). Recall that  $x_k(t)\in C^k(t)$ for all $t\in[0,T]$, and 
\begin{equation}\label{xi_constant}
|\dot{x}_k(t)|\leq ({\rm const})\;\;\;{\rm and}\;\;\; \xi_k^i(t)=\gamma_k e^{\gamma_k(\psi^i(t,x_k(t))-\sigma_k)}\leq ({\rm const}).
\end{equation}
Then  there exist subsequences (we do not relabel) weakly-$*$ converging in $L^{\infty}$ to some  $v$ and $\xi^i$.
Hence
$$x_{k}(t)=x_0+  \int_0^t \dot x_{k}(s) ds \longrightarrow x(t)=x_0+
  \int_0^t  v(s)ds, ~\forall ~t\in [0,T],
$$
for an absolutely continuous function $x$. Obviously, $x(t)\in C(t)$ for all $t\in [0,T]$.
Considering the sequence $\{x_k\}$,  recall that 
\begin{equation}
\dot{x}_k(t)\in f(t,x_k(t),U)-\sum_{i=1}^I\xi_k^i(t)\nabla_x \psi^i(t,x_k(t)).\label{ap1}
\end{equation}
Inclusion (\ref{ap1}) is equivalent to 
$$
\langle z, \dot{x}_k(t)\rangle\leq S(z,f(t,x_k(t),U))-\sum_{i=1}^I\xi_k^i(t)\langle z,
 \nabla_x\psi^i(t,x_k(t))\rangle,\;\;\;\forall\: z\in R^n.
 $$
Integrating this inequality, we get
\begin{multline}
 \left\langle z,\frac{x_k(t+\tau)-x_k(t)}{\tau}\right\rangle \\
 \leq \frac{1}{\tau}\int_t^{t+\tau}\left(S(z,f(s,x_k(s),U))-\sum_{i=1}^I\xi_k^i(s)\langle z,
 \nabla_x\psi^i(s,x_k(s))\rangle \right)ds\\
 =\frac{1}{\tau}\int_t^{t+\tau}\left(S(z,f(s,x_k(s),U))-\sum_{i=1}^I\xi_k^i(s)\langle z, \nabla_x\psi^i(s,x(s))\rangle \right.
\\
 +
\left.\sum_{i=1}^I\xi_k^i(s)\langle z, \nabla_x\psi^i(s,x(s))- \nabla_x\psi^i(s,x_k(s))\rangle \right)ds.
\end{multline}
 Passing to the limit as $k\rightarrow\infty$, we obtain
\begin{multline}
 \left\langle z,\frac{x(t+\tau)-x(t)}{\tau}\right\rangle\\ \leq \frac{1}{\tau}\int_t^{t+\tau}\left(S(z,f(s,x(s),U))-\sum_{i=1}^I\xi^i(s)\langle z,
 \nabla_x\psi^i(s,x(s))\rangle \right)ds.
\end{multline}
Let $t\in[0,T]$ be a Lebesgue point of  $x$ and $\xi$. Passing in the last inequality to the limit as $\tau\downarrow 0$, it leads to 
$$
\langle z,\dot{x}(t)\rangle\leq S(z,f(t,x(t),U))-\sum_{i=1}^I\xi^i(t)\langle z,
 \nabla_x\psi^i(t,x(t))\rangle.
 $$
Since $z\in R^n$ is an arbitrary vector and the set $f(t,x(t),U)$ is convex, we conclude that
$$
\dot{x}(t)\in f(t,x(t),U)-\sum_{i=1}^I\xi^i(t)\nabla_x\psi^i(t,x(t)).
$$
By the Filippov lemma there exists a measurable control $u(t)\in U$ such that
$$
\dot{x}(t)= f(t,x(t),u(t))-\sum_{i=1}^I\xi^i(t)\nabla_x\psi^i(t,x(t)).
$$
Furthermore, observe   that $\xi^i$ is zero if $\psi^i(t,x(t))<0$.
If for some $u$ such that $u(t)\in U$ a.e.,  $u_k=u$  for all $k$,  then the sequence $x_k$ converges to the solution of
$$
\dot{x}(t)= f(t,x(t),u(t))-\sum_{i=1}^I\xi^i(t)\nabla_x\psi^i(t,x(t)).
$$
Indeed, to see this, it suffices to pass to the limit as $k\rightarrow\infty$ and then as $\tau\downarrow 0$,  in the equality
$$
\frac{x_k(t+\tau)-x_k(t)}{\tau}= \frac{1}{\tau}\int_t^{t+\tau}\left(f(s,x_k(s),u(s))-\sum_{i=1}^I\xi_k^i(s) \nabla_x\psi^i(s,x_k(s))\right)ds.
$$

\medskip

We now prove the uniqueness of the solution. We follow the proof of Theorem 4.1 in \cite{zeidan2020}. Notice, however, that we now consider a special case and not the general case treated  in \cite{zeidan2020}.  Suppose that there exist two different solutions of (\ref{e2}): $x_1$ and $x_2$. We have
\begin{multline}
\frac{1}{2}\frac{d}{dt}|x_1(t)-x_2(t)|^2=\langle x_1(t)-x_2(t),\dot{x}_1(t)-\dot{x}_2(t)\rangle
\\
=\langle x_1(t)-x_2(t),f(t,x_1(t),u(t))-f(t,x_2(t),u(t))\rangle
\\
-\left\langle x_1(t)-x_2(t),\sum_{i=1}^I\xi_1^i(t)\nabla\psi^i(t,x_1(t))-\sum_{i=1}^I\xi_2^i(t)\nabla\psi^i(t,x_2(t))\right\rangle.
\end{multline}
If, for all $i$, $\psi^i(t,x_1(t))<0$ and $\psi^i(t,x_2(t))<0$, then $\xi_1^i(t)=\xi_2^i(t)=0$ and we obtain
$$
\frac{1}{2}\frac{d}{dt}|x_1(t)-x_2(t)|^2\leq L_f|x_1(t)-x_2(t)|^2.
$$
Suppose that $\psi^j(t,x_1(t))=0$. Then by the Taylor formula we get
\begin{multline}
\psi^j(t,x_2(t))=\psi^j(t,x_1(t))+\langle \nabla_x\psi^j(t,x_1(t)),x_2(t)-x_1(t)\rangle
\\
+\frac{1}{2}\langle x_2(t)-x_1(t), \nabla_x^2\psi^j(t,\theta x_2(t)+(1-\theta )x_1(t))( x_2(t)-x_1(t))\rangle,
\end{multline}
where $\theta\in [0,1]$. Since $\psi^j(t,x_2(t))\leq  0$, we have
\begin{multline}
\langle \nabla_x\psi^j(t,x_1(t)),x_2(t)-x_1(t)\rangle
\\
\leq -
\frac{1}{2}\langle x_2(t)-x_1(t), \nabla_x^2\psi^j(t,\theta x_2(t)+(1-\theta )x_1(t))( x_2(t)-x_1(t))\rangle\\
\leq ({\rm const})
|x_1(t)-x_2(t)|^2.
\end{multline}
Now, if $\psi^j(t,x_2(t))=0$, we deduce in the same way that
$$
\langle \nabla_x\psi^j(t,x_2(t)),x_1(t)-x_2(t)\rangle\leq ({\rm const})|x_1(t)-x_2(t)|^2.
$$
Thus we have
$$
\frac{1}{2}\frac{d}{dt}|x_1(t)-x_2(t)|^2\leq ({\rm const})|x_1(t)-x_2(t)|^2.
$$
Hence $|x_1(t)-x_2(t)|=0$. \hfill
$\Box$

\section{Approximating Family of Optimal Control Problems}
\label{ApFam}
In this section we define an approximating family of optimal control problems to $(P)$ and we state the corresponding necessary conditions. 

Let $(\hat{x},\hat{u})$  
be a global solution to $(P)$ and consider  sequences $\{\gamma_k\}$ and $\{\sigma_k\}$ as defined above.
Let  $\hat  x_k(\cdot)$ be the solution to
\begin{equation}\label{starstar}\left\{\begin{array}{rcl}
&&\dot{x}(t) =  f(t,x(t),\hat u(t))-\displaystyle\sum_{i=1}^I\gamma_k e^{\gamma_k(\psi^i(t,x(t))-\sigma_k)}\nabla_x\psi^i(t,x(t)),\\
&&x(0)=  \hat x(0).
\end{array}
\right.
\end{equation}
Set   $\epsilon_k=|\hat x_k(T)-\hat x(T)|$. 
It follows from  Theorem \ref{T1}  that $\epsilon_k\downarrow 0$. 
Take $\alpha>0$ and define the problem
$$
(P_k^\alpha)
\left\{
\begin{array}{l}
\mbox{Minimize } \; 
\phi(x(T))+|x(0)-\hat{x}(0)|^2+\alpha  \displaystyle\int_0^T|u(t)-\hat{u}(t)|dt\\[2mm]
\mbox{over processes $(x,u)$ such that }\\[2mm]
\hspace{8mm} \dot x(t) = f(t,x(t),u(t))-\displaystyle\sum_{i=1}^I\nabla_x e^{\gamma_k(\psi^i(t,x(t))-\sigma_k)} \hspace{0.2cm}\mbox{a.e.}\ \ t\in[0,T],\\[2mm]
\hspace{8mm} u(t)\in U~~~~\mbox{a.e.}\ \ t\in[0,T],\\[2mm]
\hspace{8mm} x(0)\in C_0,~~ x(T)\in C_T+\epsilon_k B_n,
\hspace{2mm}\end{array}
\right.
$$
Clearly, the problem 
$(P_k^\alpha)$ has admissible solutions.
Consider  the space 
$$
W=\{ (c,u)\mid c\in C_0,\; u\in L^{\infty}\;   {\rm with }\; u(t)\in U\}
$$
and  the distance
$$
d_{W}((c_1,u_1),(c_2,u_2))=|c_1-c_2|+\int_0^T|u_1(t)-u_2(t)|dt.
$$
Endowed with $d_{W}$,  $W$ is a complete metric space.  Take any $(c,u)\in W$ and a solution $y$ to the Cauchy problem
$$
\left\{\begin{array}{rcl}
 \dot y(t) & = &  f(t,y(t),u(t))-\displaystyle\sum_{i=1}^I\nabla_x e^{\gamma_k(\psi^i(t,y(t))-\sigma_k)} \hspace{0.2cm}\mbox{a.e.}\ \ t\in[0,T],\\[2mm]
 y(0) & = & c.
 \end{array}
 \right.
 $$
Under our assumptions,
the function 
$$
(c,u)~\rightarrow ~  \phi(y(T))+ | c - \hat{x}(0) |^2+\alpha \int_{0}^{T} 
| u-\hat{u}|~dt
$$
  is continuous on $(W,d_{W})$ and bounded below. 
 Appealing to Ekeland's Theorem we deduce the existence of a pair $(x_k,u_k)$ solving the following problem
$$
(AP_k)
\left\{
\begin{array}{l}
\mbox{Minimize } \; 
\Phi(x,{u})= \phi(x(T))+|x(0)-\hat{x}(0)|^2+\alpha  \displaystyle\int_0^T|u(t)-\hat{u}(t)|dt\\[3mm]
\qquad\qquad +\epsilon_k\left(|x(0)-x_k(0)|+  \displaystyle\int_0^T|u(t)-u_k(t)|dt\right),\\[2mm]
\mbox{over processes $(x,u)$ such that }\\[2mm]
\hspace{8mm} \dot x(t) = f(t,x(t),u(t))-\displaystyle\sum_{i=1}^I\nabla_x e^{\gamma_k(\psi^i(t,x(t))-\sigma_k)} \hspace{0.2cm}\mbox{a.e.}\ \ t\in[0,T],\\[2mm]
\hspace{8mm} u(t)\in U~~~~\mbox{a.e.}\ \ t\in[0,T],\\[2mm]
\hspace{8mm} x(0)\in C_0,~~ x(T)\in C_T+\epsilon_k B_n,
\hspace{2mm}\end{array}
\right.
$$
\begin{lemma}\label{lemma:approx}
Take $\gamma_k\to \infty$, {\color{blue}$\sigma_k\to 0$}  and  $\epsilon_k \to 0$ as defined above. For each $k$, let  $(x_k,u_k)$ be the solution to $(AP_k)$. Then there exists a  subsequence (we do not relabel) such that 
$$
u_k(t)\rightarrow \hat u(t)~  {a.e.}, \quad x_k\rightarrow \hat x\;  \rm{ uniformly}\; in\;  [0,T].
$$
\end{lemma}
\textit{ Proof}
We deduce  from Theorem \ref{T1} that $\{x_k\}$ uniformly converges to an admissible  solution $\tilde{x}$ to  
$(P)$.  
Since $U$ and $C_0$ are compact, we have $U\subset KB_m$ and $C_0\subset KB_n$. Without loss of generality, $u_k$ weakly-$*$ converges to a function $\tilde{u}\in L_{\infty}([0,T],U)$. Hence it weakly  converges  to $\tilde{u}$ in $L_1$.  From optimality of the processes $(x_k,u_k)$ we have
$$
\phi(x_k(T))+|x_k(0)-\hat{x}(0)|^2+\alpha\int_0^T|u_k(t)-\hat{u}(t)|dt
$$
$$
\leq \phi(\hat{x}_k(T))+\epsilon_k\left(
|\hat{x}_k(0)-x_k(0)|+\int_0^T|u_k(t)-\hat{u}(t)|dt\right)
$$
$$
\leq \phi(\hat{x}_k(T))+2K(1+T)\epsilon_k.
$$
 Since $(\hat{x},\hat{u})$ is  a global solution of the problem, passing to the limit, we get
$$
\phi(\tilde{x}(T))+|\tilde{x}(0)-\hat{x}(0)|^2+\alpha\int_0^T|\tilde{u}(t)-\hat{u}(t)|dt
$$
$$
\leq\lim_{k\rightarrow\infty}(\phi(x_k(T))+|x_k(0)-\hat{x}(0)|^2)+ \alpha\liminf_{k\rightarrow\infty}
\int_0^T|u_k(t)-\hat{u}(t)|dt
$$
$$
\leq \lim_{k\rightarrow\infty}\phi(\hat{x}_k(T))=\phi(\hat{x}(T))\leq \phi(\tilde{x}(T)).
$$
Hence $\tilde{x}(0)=\hat{x}(0)$, $\tilde{u}=\hat{u}$ a.e., and $u_k$ converges to $\hat{u}$ in $L_1$, and some subsequence converges to $\hat{u}$ almost everywhere (we do not relabel).\hfill $\Box$

\bigskip

We now finish this section with  the statement of the optimality necessary conditions for the family of problems $(AP_k)$. These can be seen as a  direct consequence of  Theorem 6.2.1 in \cite{Vinter}.
\begin{proposition}\label{lemma:nco}
For each k, let $(x_k,u_k)$ be a solution to $(AP_k)$. Then there exist absolutely continous functions $p_k$ and scalars  $\lambda_k\geq 0$ such that
\begin{itemize}
\item[{\bf(a)}] (nontriviality condition) \begin{equation}
\label{not0}
  \lambda_k+|p_{k}(T)|  =1,
\end{equation}

\item[{\bf(b)}] (adjoint equation)
\begin{equation}
\label{dynamic}
\begin{array}{c}\dot{p}_{k} = {\color{red}-}(\nabla_x f_{k})^* p_{k} 
+\sum_{i=1}^I\gamma_k e^{\gamma_k (\psi_{k}^i-\sigma_k)}\nabla^2_x\psi_{k}^ip_{k}\\[2mm]
+\sum_{i=1}^I\gamma_k^2e^{\gamma_k (\psi_{k}^i-\sigma_k)}\nabla_x\psi_{k}^i\langle\nabla_x\psi_{k}^i,p_{k}\rangle, 
\end{array}
\end{equation}
where the superscript $*$ stands for transpose,

\item [{\bf(c)}] (maximization condition)
\begin{equation}
\label{max}
   \max_{u\in U}\left\{ \langle f(t,x_{k}, u)  ,  p_{k} \rangle -  \alpha \lambda_k|u-\hat{u}|
-\epsilon_k \lambda_k|u-u_k|\right\}
\end{equation}
is attained at  $u_k (t) $, for almost every $t\in [0,T]$,
\item [{\bf(d)}] (transversality condition)
\begin{eqnarray}
 ( p_{k}(0), - p_{k}(T))
 \in \lambda_k\left(2(x_k(0)-\hat{x}(0))+\epsilon_k B_n,  \partial  \phi (x_{k}(T))\right)\nonumber \\[2mm] 
 +  N_{C_0}(x_{k}(0))\times N_{C_T+\epsilon_kB_n}(x_{k}(T)). \hspace{1cm}\label{trans}
\end{eqnarray}
%\item [{\bf(iv)}] (nontriviality condition)
%\begin{equation}
%\label{not0}
% |p_{\gamma_k}(1)| + \lambda_k =1.
%\end{equation}
\end{itemize}
\end{proposition}

To simplify the notation above, we drop the $t$ dependance in $p_k$, $\dot p_k$, $x_k$, $u_k$, $\hat x$ and $\hat u$. Moreover,  in (b), we write $\psi_k$ instead of $\psi(t,x_k(t))$, $f_k$ instead of $f(t,x_k(t),u_k(t))$. The same holds for the derivatives of $\psi$ and $f$.

\section{Maximum Principle for $(P)$}
\label{MP}
In this section, we establish our main result, a Maximum Principle for $(P)$.  
This is done by taking  limits of the conclusions of Proposition \ref{lemma:nco}, following closely the analysis done in the proof of \cite[Theorem 2]{nosso_2019}.

Observe  that
\begin{equation}\nonumber
\begin{split}
\frac{1}{2} \frac{d}{dt} |p_k(t)|^2  & = 
- \langle \nabla_x f_k p_k , p_k \rangle +  
\sum_{i=1}^I\gamma_k e^{\gamma_k(\psi_k^i-\sigma_k)} \langle \nabla_x^2\psi_k^ip_k, p_k \rangle \\ 
& \hspace{1cm} +  
\sum_{i=1}^I\gamma_k^2e^{\gamma_k(\psi_k^i-\sigma_k)}\langle \nabla_x\psi_k^i, p_k \rangle ^2 
\\
& \geq - \langle \nabla_x f_k p_k , p_k \rangle + \sum_{i=1}^I 
\gamma_k e^{\gamma_k(\psi_k^i-\sigma_k)} \langle \nabla_x^2\psi_k^ip_k, p_k \rangle
\\
& \geq \ - M | p_k|^2 +\sum_{i=1}^I\gamma_k e^{\gamma_k(\psi_k^i-\sigma_k)} \langle \nabla_x^2\psi_k^ip_k, p_k \rangle,
\end{split}
\end{equation}
where $M$ is the constant of (A2). Taking into account hypothesis  (A1) and (\ref{star}) we deduce the existence of a constant $K_0>0$ such that
$$
\frac{1}{2} \frac{d}{dt} |p_k(t)|^2\geq  -K_0| p_k{\color{blue}(t)}|^2.
$$
This last inequality leads to 
$$
| p_k(t)|^2 \ \leq \ e^{2 K_0 (T-t)} | p_k(T)|^2 \leq  \ e^{2 K_0T} |p_k(T)|^2.
$$
Since, by (a) of Proposition \ref{lemma:nco},  $|p_k(T)|\leq 1$, we deduce from the above that
there exists  $M_0>0$ such that
\begin{equation}\label{pk1}
| p_k(t)| \ \leq  M_0.
\end{equation}
Now, we claim that  the sequence $\{\dot p_k\}$ is uniformly bounded in $L^1$.  To prove  our claim, we need to establish  bounds  for the three terms in \eqref{dynamic}. Following  \cite{nosso_2019} and  \cite{nosso_2022}, we start  by   deducing  some inequalities that  will be of help.

Denote $I_k=I(t,x_k(t))$ and  $S_k^j={\rm sign}\left( \langle  \nabla_x\psi_k^j, p_ k \rangle \right)$. We have
\begin{equation}\nonumber
\begin{split}
\sum_{j=1}^I\frac{d}{dt} & \left| \langle    \nabla_x\psi_k^j, p_ k \rangle \right|  \\
& =\sum_{j=1}^I\left( \langle \nabla^2_x\psi_k^j \dot{x}_ k, p_ k \rangle
+\langle\partial_t\nabla_x\psi_k^j,p_k\rangle
+ \langle\nabla_x\psi_k^j,\dot{p}_ k\rangle \right)\, S_k^j
% {\rm sign}\left( \langle  \nabla_x\psi_k^j, p_ k \rangle \right)
\\
& =
\sum_{j=1}^I\left( \langle  p_ k, \nabla^2_x\psi_k^j f_ k \rangle
-  \sum_{i=1}^I\gamma_k e^{\gamma_k(\psi_k^i-\sigma_k)}  \langle p_ k, \nabla^2\psi_k^j   \nabla_x\psi_k^i  \rangle \right)S_k^j
%{\rm sign}\left( \langle  \nabla_x\psi_k^j, p_ k \rangle \right)
\\
& \hspace{1cm}+\sum_{j=1}^I\left( \langle\partial_t\nabla_x\psi_k^j,p_k\rangle
- \langle  \nabla_x \psi_k^j, (\nabla_x f_ k)^* p_ k \rangle \right)S_k^j
%{\rm sign}\left( \langle  \nabla_x\psi_k^j, p_ k \rangle \right)
\\ 
& \hspace{1cm}+\sum_{j=1}^I\left( 
\sum_{i=1}^I\gamma_k e^{\gamma_k(\psi_k^i-\sigma_k)}  \langle  \nabla_x\psi_k^j ,  \nabla^2_x\psi_k^i  p_ k\rangle\right)S_k^j
%{\rm sign}\left( \langle  \nabla_x\psi_k^j, p_ k \rangle \right)
\\
& \hspace{1cm}+\sum_{i=1}^I \sum_{j=1}^I \gamma_k^2e^{\gamma_k(\psi_k^i-\sigma_k)}  \langle  \nabla_x\psi_k^j ,  \nabla_x\psi_k^i \rangle \langle\nabla_x\psi_k^i,p_ k\rangle  S_k^j
% {\rm sign}\left( \langle  \nabla_x\psi_k^j, p_ k \rangle \right).
\end{split}
\end{equation}
Observe that (see (\ref{cn2}) and \eqref{gradient0})
\begin{multline}\nonumber
\sum_{i=1}^I \sum_{j=1}^I \gamma_k^2e^{\gamma_k(\psi_k^i-\sigma_k)}  \langle  \nabla_x\psi_k^j ,  \nabla_x\psi_k^i \rangle \langle\nabla_x\psi_k^i,p_ k\rangle  S_k^j
% {\rm sign}\left( \langle  \nabla_x\psi_k^j, p_ k \rangle \right)
\\
=\sum_{i=1}^I \sum_{j\in I_k} \gamma_k^2e^{\gamma_k(\psi_k^i-\sigma_k)}  \langle  \nabla_x\psi_k^j ,  \nabla_x\psi_k^i \rangle \langle\nabla_x\psi_k^i,p_ k\rangle  S_k^j
% {\rm sign}\left( \langle  \nabla_x\psi_k^j, p_ k \rangle \right)
\\
=\sum_{i\not\in I_k} \gamma_k^2e^{\gamma_k(\psi_k^i-\sigma_k)} \sum_{j\in I_k}  \langle  \nabla_x\psi_k^j ,  \nabla_x\psi_k^i \rangle \langle\nabla_x\psi_k^i,p_ k\rangle  
S_k^j
% {\rm sign}\left( \langle  \nabla_x\psi_k^j, p_ k \rangle \right)
\\
+\sum_{i\in I_k}\gamma_k^2e^{\gamma_k(\psi_k^i-\sigma_k)}\left(
|\nabla_x\psi_k^i|^2+
 \sum_{j\in I_k\setminus\{i\}}   \langle  \nabla_x\psi_k^j ,  \nabla_x\psi_k^i \rangle
 S_k^j~S_k^i
%  \, {\rm sign}\left( \langle  \nabla_x\psi_k^j, p_ k \rangle \right)
 %  \, {\rm sign}\left( \langle  \nabla_x\psi_k^i, p_ k \rangle \right)
   \right)
  | \langle\nabla_x\psi_k^i,p_ k\rangle|
\\
=  \sum_{i\in I_k}\gamma_k^2e^{\gamma_k(\psi_k^i-\sigma_k)}\left(
|\nabla_x\psi_k^i|^2+
 \sum_{j\in I_k\setminus\{i\}}   \langle  \nabla_x\psi_k^j ,  \nabla_x\psi_k^i \rangle
 S_k^j~S_k^i
%  \, {\rm sign}\left( \langle  \nabla_x\psi_k^j, p_ k \rangle \right)
 %  \, {\rm sign}\left( \langle  \nabla_x\psi_k^i, p_ k \rangle \right)
   \right)
  | \langle\nabla_x\psi_k^i,p_ k\rangle|\\
\geq \disp(1-\rho)\sum_{i\in I_k}\gamma_k^2e^{\gamma_k(\psi_k^i-\sigma_k)}|\nabla_x\psi_k^i|^2 | \langle\nabla_x\psi_k^i,p_ k\rangle| 
\\
= 
(1-\rho)\sum_{i=1}^I\gamma_k^2e^{\gamma_k( \psi_k^i-\sigma_k)}|\nabla_x\psi_k^i|^2 | \langle\nabla_x\psi_k^i,p_ k\rangle|. \hspace{2cm}
\end{multline}
%The last equality above is due to \eqref{gradient0}.

Using this and integrating the previous equality, we deduce the existence of $M_1>0$ such that:
\begin{equation}\label{bound:int}
\int_0^T \sum_{i=1}^I\gamma_k^2e^{\gamma_k(\psi_k^i-\sigma_k)}|  \nabla_x\psi_k^i|^2 | \langle\nabla_x\psi_k^i,p_ k\rangle|dt
 \leq M_1.
\end{equation}

We are now in a position to show   that  $$  \disp\int_{0}^{T} \sum_{i=1}^{I} \gamma_k^2  e^{\gamma_k( \psi_k^i-\sigma_k)} |\nabla_x \psi_k^i|  \left| \langle \nabla_x \psi_k^i,p_ k\rangle \right| \ dt $$ is bounded.
For simplicity, set 
$L_k^i(t) =\gamma_k^2 e^{\gamma_k( \psi_k^i-\sigma_k)} |\nabla_x \psi_k^i|  \left| \langle \nabla_x \psi_k^i,p_ k\rangle\right|$. Notice  that
$$  \disp\sum_{i=1}^I \int_{0}^{T}   L_k^i(t) dt=   \disp\sum_{i=1}^I \left\lbrace \int_{\{t: |\nabla_x \psi_k^i| <  \eta\}} \hspace{-0.7cm}L_k^i(t)~dt+
\disp  \int_{\{t: |\nabla_x\psi_k^i| \geq \eta\}}\hspace{-0.7cm} L_k^i(t) dt \right\rbrace.$$
Using (A1)  and \eqref{bound:int}, we deduce that 
\begin{equation}\nonumber
\begin{split}
  \disp\sum_{i=1}^I  \int_{0}^{T} L_k^i (t)~dt  
& \leq    \disp\sum_{i=1}^I  \left (\gamma_k^2  e^{-\gamma_k(\beta+\sigma_k)} \eta^2 \max_{t} |p_ k(t)|\right)
\\
& \hspace{0.5cm}+\disp\sum_{i=1}^I \left(
\gamma_k^2  \int_{ \{t: |\nabla_x \psi_k^i | \geq \eta\} } \hspace{-1cm}e^{\gamma_k(\psi_k^i-\sigma_k)} \frac{|\nabla_x\psi_k^i |^2}{ | \nabla_x\psi_k^i |} \left| \langle  \nabla_x\psi_k^i,p_ k\rangle \right| \ dt \right)
\\[1mm]
& \leq \gamma_k^2  I~e^{-\gamma_k (\beta+\sigma_k)} \eta^2 M_0  \\
&\hspace{0.5cm} +\frac{1}{\eta}\disp\sum_{i=1}^I \left(
  \int_{0}^{T} \gamma_k^2 e^{\gamma_k(\psi_k^i-\sigma_k)} | \nabla_x\psi_k^i|^2 \left| \langle \nabla_x\psi_k^i,p_ k\rangle \right|\ dt\right)
\\[1mm]
& \leq \ \eta^2 M_0 {\color{blue}I}+ \frac{M_1}{\eta},
\end{split}
\end{equation}
for $k$ large  enough. 
Summarizing, there exists a  $M_2>0$ such that 
\begin{equation}
\label{M2}
\disp\sum_{i=1}^I  \gamma_k^2 \int_{0}^{T} e^{\gamma_k(\psi_k^i-\sigma_k)} |\nabla\psi_k^i|  \left| \langle \nabla\psi_k^i,p_ k\rangle \right| \ dt \  \ \leq  M_2.
\end{equation}
Mimicking the analysis conducted in Step 1, b) and c) of  the proof of Theorem 2 in \cite{nosso_2019}
 and taking into account (b) of Proposition \ref{lemma:nco} we conclude that there exist constants $N_1>0$ %and $N_2>0$ 
 such that

\begin{equation}
\label{M7}
\int_{0}^{T} \left| \dot{p}_{\gamma_k}(t)\right|dt \leq N_1, 
\end{equation}
for $k$ sufficiently large, proving our claim.

Before proceeding, observe that it is a simple matter to assert the existence of a constant $N_2$
such that
 \begin{equation}
\label{M8}
\disp\sum_{i=1}^I  \int_0^T \gamma_k^2 e^{\gamma_k(\psi_k^i-\sigma_k)} |\langle\nabla\psi_k^i,p_{\gamma_k}\rangle|dt\leq N_2.
\end{equation}
This inequality will be of help in what follows.

 Let us now recall  that $$\xi_k^i(t)=\gamma_k e^{\gamma_k( \psi^i(t,x_k(t))-\sigma_k)}$$ 
 and that the second inequality in \eqref{xi_constant} holds.
We turn to  the  analysis of Step 2  in the proof of Theorem 2 in \cite{nosso_2019} (see also \cite{nosso_2022}). Adapting those arguments,  we can   conclude the existence of some 
  function  $p\in BV([0,T],R^n)$ and, for $i=1, \ldots, I$, functions   $\xi^i\in L^{\infty}([0,T],R)$ with  $\xi^i(t) \geq 0 \ \mbox{ a. e. } t$, $\xi^i(t) = 0, \ t \in I_b^i$,  where
$$I_b^i=\left\{t\in [0,T]:~ \psi^i(t, \hat x(t))<0\right\},$$
 and  finite signed Radon measures $\eta^i$, null in $ I_b^i$, such that, for any $z\in C([0,T],R^n)$
$$
\int_0^T \langle z,dp\rangle=-\int_0^T \langle z, (\nabla\hat{f})^*p\rangle dt +\disp\sum_{i=1}^{I} \left( \int_0^T\xi^i\langle z,\nabla^2\hat{\psi} ^ip\rangle dt
+\int_0^T\langle z, \nabla\hat{\psi}^i (t)  \rangle d\eta^i\right), 
$$
where $\nabla \hat \psi^i(t)=\nabla   \psi^i(t,\hat x(t))$. 
The finite signed Radon measures $\eta^i$ are weak-$*$ limits of 
$$\gamma_k^2 e^{\gamma_k(\psi_k^i-\sigma_k)} \langle\nabla\psi_k^i(x_k(t),p_{k}(t)\rangle dt.$$ Observe that the measures 
\begin{equation}\label{mea:nova}
\langle \nabla\psi^i(\hat x(t),p(t)\rangle d\eta^i(t)
\end{equation}
are nonnegative.

For each $i=1,\ldots, I$, the  sequence $\xi_k^i$ is weakly-$*$ convergent  in $L^{\infty}$ to $\xi^i\geq 0$. Following \cite{nosso_2022}, we deduce from  (\ref{M8}) that, for each $i=1,\ldots, I$, 
$$
\int_0^T|\xi^i\langle\nabla_x\hat{\psi}^i ,p\rangle|  dt=\lim_{k\rightarrow\infty}\int_0^T|\xi_k^i\langle\nabla_x\hat{\psi}^i,p\rangle| dt
$$
$$
\leq \lim_{k\rightarrow\infty} \left( \int_0^T\xi_k^i|\langle\nabla_x\hat{\psi}^i ,p\rangle-\langle\nabla_x \psi_k^i,p_k\rangle| dt+
\int_0^T\xi_k^i|\langle\nabla_x \psi_k^i,p_k\rangle| dt\right)
$$
$$
\leq \lim_{k\rightarrow\infty}\left( \Big|\xi_k^i\Big|_{L^{\infty}}\Big| \langle\nabla_x\hat{\psi}^i,p\rangle-\langle \nabla_x \psi_k^i,p_k\rangle\Big|_{L^1}+\frac{N_2}{\gamma_k}\right)=0.
$$
It turns out  that 
\begin{equation}
\label{ccc}
\xi^i \langle\nabla_x\hat{\psi}^i ,p\rangle=0\; {\rm a. e.}.
\end{equation}

Consider now the sequence of scalars $\{\lambda_k\}$. It is an easy matter to show  that there exists a subsequence of  $\{\lambda_k\}$ converging to some $\lambda\geq 0$. This, together with the  convergence of $p_k$ to $p$,  allows us to take limits in (a) and (c) of Proposition \ref{lemma:nco} to deduce that
$$\lambda+|p(T)|=1$$
and
$$\langle p(t), f(t,\hat x(t),u)\rangle-\alpha\lambda |u-\hat{u}(t)| \leq \langle p(t), f(t,\hat x(t),\hat u(t))\rangle  ~\forall u\in U,  \text{ a.e. } t \in [0,T].$$

It remains to take limits of the transversality conditions (d) in Proposition \ref{lemma:nco}.
First, observe that 
$$C_T+\epsilon_kB_n=\left\{x:~d(x,C_T)\leq \epsilon_k\right\}.$$
From the basic properties of the Mordukhovich normal cone and subdifferential (see \cite{Mor05}, section 1.3.3) we have
$$
N_{C_T+\epsilon_kB_n}(x_k(T))\subset   \text{ cl cone}\:\partial d(x_k(T), C_T)$$
and
$$
N_{C_T}(\hat{x}(T))=   \text{ cl cone}\:\partial d(\hat{x}(T), C_T).
$$
Passing to the limit as $k\to\infty$ we get
$$(p(0),-p(T))\in 
N_{C_0}(\hat x(0))\times N_{C_T}(\hat x(T))+\{0\}\times \lambda ~ \partial  \phi(\hat x(T)).$$

Finally, and mimicking  Step 3 in the proof of Theorem 2 in \cite{nosso_2019}, we  remove the dependence of the conditions on the parameter $\alpha$. This is done by taking further limits, this time  considering a sequence of $\alpha_j\downarrow 0$. 

 We then summarize our conclusions in the  following Theorem.
%
%{\color{red} The Lemma is written here but I cannot explain clearly  how this helps to take limits of the tranversality conditions.}
%\begin{lemma}
%Let $x_0\in C_0$ and $\psi(0,x_0)=0$. Assume that $-\nabla_x\psi (0,x_0)\not\in N^L_{C_0}(x_0)$. Let
%$q_k\in N^L_{(C_0+\epsilon_k B_n)\cap C^k(0)}(x_k)$, $\epsilon_k\downarrow 0$, $x_k\rightarrow x_0$, $q_k\rightarrow q$.
%Then $q\in N^L_{C_0}(x_0)+R_+\nabla_x\psi (0,x_0)$.
%\end{lemma}
%
%{\em Proof}. Show that $-\nabla_x\psi (0,x_k)\not\in N^L_{(C_0+\epsilon_k B_n)}(x_k)$. Suppose that
%$$
% -\nabla_x\psi (0,x_k)\in N^L_{(C_0+\epsilon_k B_n)\cap C^k(0)}(x_k). 
% $$
% Since
%$$
%  N^L_{(C_0+\epsilon_k B_n)}(x_k)\subset {\rm cl\: cone\:}\partial d(x_k,C_0),
%  $$
%  we have
% $$
%-\nabla_x\psi(0,x_0)=-\lim_{k\rightarrow\infty} \nabla_x\psi(0,x_k)\in 
%\limsup_{k\rightarrow\infty}{\rm cl\: cone\:}\partial d(x_k,C_0)\subset {\rm cone\:} \partial d(x_0,C_0)=N^L_{C_0}(x_0),
% $$
% a contradiction. Hence
% $$
% q_k\in  N^L_{(C_0+\epsilon_k B_n)\cap C^k(0)}(x_k)\subset N^L_{C_0+\epsilon_k B_n}(x_k)+R_+\nabla_x\psi (0,x_k)
% ={\rm cl\: cone\:}\partial d(x_k,C_0)+R_+\nabla_x\psi (0,x_k).
% $$
% Therefore, since $N^L_{C^k(0)}(x_k)=R_+\nabla_x\psi (0,x_k)$, we get
% $$
% q=\lim_{k\rightarrow\infty} q_k \in \limsup_{k\rightarrow\infty}({\rm cl\: cone\:}\partial d(x_k,C_0)+R_+\nabla_x\psi (0,x_k))
% $$
% $$
% \subset
% {\rm cone\:} \partial d(x_0,C_0)+R_+\nabla_x\psi (0,x_0)=N^L_{C_0}(x_0)+R_+\nabla_x\psi (0,x_0).\;\;\;$\Box$
% $$
% 
%\vspace{10mm}
%
%If $x_k\rightarrow x_0$, $\psi(0,x_0)<0$, then $q\in N^L_{C_0}(x_0)$.
%

\begin{theorem}\label{main}
 Let $(\hat{x}, \hat{u})$  be the optimal solution to  $(P)$. 
Suppose that assumption A1--A6 are satisfied. For $i=1,\cdots , I$, set 
$$I^{i}_b= \{ t \in [0,T]: ~ \psi^{i}(t,\hat{x}(t) ) < 0 \}.$$
There exist $ \lambda \geq 0$, $ p\in BV([0,T],R^n)$, finite signed Randon measures 
 $ \eta^i$,  null in $I^{i}_b$, for $i=1,\cdots , I$, 
$ \xi^{i}\in L^\infty([0,T],R)$, with $i=1,\cdots , I  $, where $ \disp  \xi^{i}(t) \geq 0 \ \text{ a. e. } t$ and $\xi^{i}(t) = 0, \ t \in I^{i}_b,  $ 
such that 
 
\begin{itemize}
\item[\textbf{a)}] $\lambda+|p(T)|\neq 0$,\\[1mm]
\item[\textbf{b)}] $\dot{ \hat x}(t)=f(t,\hat x(t),\hat u(t))- \disp\sum_{i=1}^{I}\xi^i(t)\nabla_x \hat \psi^{i} (t),$
\item[\textbf{c)}] for any $z\in C([0,T];R^n)$ 
$$\begin{array}{c}
\disp\int_0^T \langle z(t),dp(t)\rangle =  -\disp\int_0^T  \langle z(t), ( \nabla_x \hat f(t))^*p(t)\rangle dt
%\\[1mm]
% \hspace{-1mm} +
\\[3mm]
% \hspace{-1mm} 
% \hspace{10mm} 
 \disp  + \sum_{i=1}^{I}\disp \left( \int_0^T \xi^{i}(t) \langle z(t), \nabla^2_x\hat \psi^{i}(t) p(t)\rangle dt \right.\
 +\disp \left. \int_0^T  \langle z(t), \nabla_x \hat \psi^{i}(t)\rangle d\eta_i\right),\end{array}$$
where   $ \nabla \hat f(t) = \nabla_x f(t,\hat x(t),\hat u(t)), ~~\nabla \hat \psi^i(t)=\nabla \psi^i(t,\hat x(t))$ and $\nabla^2 \hat \psi^i(t)=\nabla^2 \psi^i(t, x(t)),$\\[0.5mm]
\item[\textbf{d)}] $\xi_i(t)\langle \nabla_x \psi^{i}(t,\hat x(t)),p(t)\rangle =0$,  $a.e. \, t$ for all $i=1, \ldots, I$,
\\[0.5mm]
\item[\textbf{e)}] for all $i=1, \ldots, I$, the meaures $\langle \nabla\psi^i(\hat x(t),p(t)\rangle d\eta^i(t)$ are nonnegative,
\\[0.5mm]
 \item[\textbf{f)}]
$\disp \langle p(t), f(t,\hat x(t),u)\rangle \leq \langle p(t), f(t,\hat x(t),\hat u(t))\rangle$  for all $u \in U,$ $~a.e.\, t$,
\\[0.5mm]
\item[\textbf{g)}]
$\disp \begin{array}{c}(p(0),-p(T))\in N_{C_0}(\hat x(0))\times N_{C_T}(\hat x(T)) 
+\{0\}\times\lambda  \partial  \phi(\hat x(T)).\end{array}$
\end{itemize}
%Then there exist a non negative scalar $\lambda$,
%$p\in BV([0,T]; R^n)$,   a finite signed Radon measure $\eta$, null in $I_b$, 
%$\xi \in L^{\infty}([0,T]; R)$ with
%$   \xi(t)\geq 0   { a.e. } t, ~\xi(t)=0   { for } t \in I_b$, 
%satisfying the following conditions
%\begin{itemize}
%\item [(a)] 
%$\lambda+|p(T)|\neq 0$,
%
%\medskip
%
%\item [(b)] 
%for any $z\in C([0,1]; R^n)$ 
%$$\begin{array}{c}
%\displaystyle  \int_0^T \langle z(t),dp(t)\rangle =
%- \displaystyle \int_0^T \langle z(t), ( \nabla_x \hat f(t))^*p(t)\rangle dt\\[2mm] 
%\hspace{1cm} +
% \displaystyle \int_0^1\!\xi(t) \langle z(t), \nabla^2_x\hat \psi(t) p(t)\rangle dt+  \int_0^1\!\langle z(t), \nabla_x \hat \psi(t)\rangle d\eta,\end{array}$$
%where the superscript $*$ stands for transpose,  $ \nabla_x \hat f(t) = \nabla_x f(t,\hat x(t),\hat u(t))$, $\nabla_x \hat \psi(t)=\nabla_x \psi(t,\hat x(t))$ and $\nabla^2_x \hat \psi(t)=\nabla^2_x \psi(t,\hat x(t)),$
%\medskip
%
%\item [(c)]
%$\langle p(t), f(t,\hat x(t),u)\rangle \leq \langle p(t), f(t,\hat x(t),\hat u(t))\rangle$ for  a.e. $t \in [0,T]$
%and  all $u\in U,$
%\medskip
%
%\item[(d)]
%$(p(0),-p(T))\in 
%N_{C_0}(\hat x(0))\times N_{C_T}(\hat x(T))+\{0\}\times \lambda  \partial  \phi(\hat x(T))
%$
%\end{itemize}
\end{theorem}

\bigskip

Noteworthy,  condition \textbf{e)} is not  considered in any of our previous works.

\bigskip

We  now turn to  the free end point case, i. e.,  to the problem 
$$
(P_f)
\left\{
\begin{array}{l}
\mbox{Minimize } \; 
\phi(x(T))\\[2mm]
\mbox{over processes $(x,u)$ such that }\\[2mm]
\hspace{8mm} \dot x(t) \in f(t,x(t),u(t))- N_{C(t)}(x(t)), \hspace{0.2cm}\mbox{a.e.}\ \ t\in[0,T],\\[2mm]
\hspace{8mm} u(t)\in U, \ \ \;\, \mbox{a.e.}\ \ t\in[0,T],\\[2mm]
\hspace{8mm} x(0) \in  C_0 \subset C(0).
\hspace{2mm}\end{array}
\right.
$$
Problem  $(P_f)$ differs from $(P)$  because  $x(T)$  is not constrained to take values in $C_T$. 
We apply  Theorem \ref{main} to $(P_f)$. Since $x(T)$ is free,  we deduce from (f) in the above Theorem that  $-p(T)=\lambda \partial \phi(\hat x(T))$. Suppose that  $\lambda=0$. Then $p(T)=0$  contradicting the nontriviality condition  (a) of Theorem \ref{main}. Without loss of generality, we then  conclude that the conditions of  Theorem \ref{main}  hold with $\lambda=1$. We summarize our findings  in the following Corollary.
\begin{corollary}\label{main_free}
 Let $(\hat{x}, \hat{u})$  be the optimal solution to  $(P_f)$. 
Suppose that assumption A1--A6 are satisfied. For $i=1,\cdots , I$, set 
$$I^{i}_b= \{ t \in [0,T]: ~ \psi^{i}(t,\hat{x}(t) ) < 0 \}.$$
There exist   $ p\in BV([0,T],R^n)$, finite signed Randon measures 
 $ \eta_i$,  null in $I^{i}_b$, for $i=1,\cdots , I$, 
$ \xi^{i}\in L^\infty([0,T],R)$, with $i=1,\cdots , I  $, where $ \disp  \xi^{i}(t) \geq 0 \ \text{ a.e. } t$ and $\xi^{i}(t) = 0$ for $t \in I^{i}_b,  $ 
such that 
 
\begin{itemize}
\item[\textbf{a)}] $\dot{ \hat x}(t)=f(t,\hat x(t),\hat u(t))- \disp\sum_{i=1}^{I}\xi^i(t)\nabla_x \hat \psi^{i} (t),$
\item[\textbf{b)}]for any $z\in C([0,T];R^n)$ 
$$\begin{array}{c}
\disp\int_0^T \langle z(t),dp(t)\rangle =  -\disp\int_0^T  \langle z(t), ( \nabla_x \hat f(t))^*p(t)\rangle dt
%\\[1mm]
% \hspace{-1mm} +
\\[3mm]
% \hspace{-1mm} 
% \hspace{10mm} 
 \disp  + \sum_{i=1}^{I}\disp \left( \int_0^T \xi^{i}(t) \langle z(t), \nabla^2_x\hat \psi^{i}(t) p(t)\rangle dt \right.\
 +\disp \left. \int_0^T  \langle z(t), \nabla_x \hat \psi^{i}(t)\rangle d\eta_i\right),\end{array}$$
where   $ \nabla \hat f(t) = \nabla_x f(t,\hat x(t),\hat u(t)), ~~\nabla \hat \psi^i(t)=\nabla \psi^i(t,\hat x(t))$ and $\nabla^2 \hat \psi^i(t)=\nabla^2 \psi^i(t, x(t)),$\\[0.5mm]
\item[\textbf{c)}] $\xi^i(t)\langle \nabla_x \psi^{i}(t,\hat x(t)),p(t)\rangle =0$ for   $a. e.~t$ and  for all $i=1, \ldots, I$,
\\[0.5mm]
\item[\textbf{d)}] for all $i=1, \ldots, I$, the meaures $\langle \nabla\psi^i(\hat x(t),p(t)\rangle d\eta^i(t)$ are nonnegative,
\\[0.5mm]
\item[\textbf{e)}] $\disp \langle p(t), f(t,\hat x(t),u)\rangle \leq \langle p(t), f(t,\hat x(t),\hat u(t))\rangle$ for all $u \in U$, $a.e.\, t$,
\\[0.5mm]
\item[\textbf{f)}]
$\disp \begin{array}{c}(p(0),-p(T))\in N_{C_0}(\hat x(0))\times \{0\}
+\{0\}\times  \partial  \phi(\hat x(T)).\end{array}$
\end{itemize}
\end{corollary}

%Observe that the Maximum Principle for problem $(P_f)$ in  \cite{Palladino2022} is not in the normal from, i.e., the cost multiplier $\lambda$ may be zero. 
%The remarkable feature of Corollary \ref{main_free} is that it holds when $\lambda=1$. 

\section{Example} 

Let us consider the following problem 
$$\left\{
\begin{array}{l}
\mbox{Minimize } \; 
 -x(T)\\[2mm]
\mbox{over processes $((x,y,z),u)$ such that }\\[2mm]
\hspace{8mm} \begin{bmatrix}
\dot x(t)\\ \dot y(t)\\ \dot z(t)
\end{bmatrix}
\in 
\begin{bmatrix}
0 & \sigma & 0\\
0 & 0 & 0\\
0 & 0 & 0
\end{bmatrix}
\begin{bmatrix}
x\\ y\\ z
\end{bmatrix}
+\begin{bmatrix}
0\\ u\\ 0
\end{bmatrix}-N_C(x,y,z),\\[2mm]
\hspace{8mm}   u\in [-1,1],\\[1mm]
\hspace{8mm}  (x,y,z)(0)=(x_0,y_0,z_0), \\[1mm]
\hspace{8mm}  (x,y,z)(T)\in C_T, 
\end{array}
\right.$$
where
\begin{itemize}
\item[$\bullet$] $0<\sigma\ll 1$,\\

\item[$\bullet$] $
C=\{ (x,y,z)\mid x^2+y^2+(z+h)^2\leq1,\; x^2+y^2+(z-h)^2\leq1\},\;\; 2h^2<1
$,\\

\item[$\bullet$] $(x_0,y_0,z_0)\in {\rm int }C$,  with  $x_0<-\delta $, $y_0=0$ and $z_0>0$,\\
%\cap \left\{(x,y,z)\in R^3: ~x<0,~y=0, ~ z>0\right\}$, 

\item[$\bullet$] $C_T=\{ (x,y,z)\mid x\leq0,\; y\geq 0,\; \delta y-y_2x\leq\delta y_2\}\cap C$, where 
$$
\delta<\frac{y_2|x_0|}{y_1},~ {\rm with } ~y_1=\sqrt{1-x_0^2-(z_0+h)^2} \text{ and }y_2=\sqrt{1-h^2}.
$$
\end{itemize}

We choose $T>0$ small and, nonetheless,  sufficiently large to guarantee that, when  $\sigma=0$,  the system  can reach the interior of $C_T$ but not  the segment  $\{ (x,0,0)\mid x\in [-\delta,0]\}$. Since $\sigma$ and $T$ are  small, it follows that the optimal trajectory should reach $C_T$ at the face $\delta y-y_2x=\delta y_2$ of $C_T$. 

To significantly increase the value of the  $x(T)$,   the optimal  trajectory needs to  live  on the boundary of $C$ for some interval of time. Then, before reaching and after leaving the boundary of $C$,   the optimal trajectory  lives in the interior of $C$. Since $\delta$ is small, the trajectory cannot reach $C_T$ from any point of the sphere $x^2+y^2+(z+h)^2=1$ with $z>0$. This means that, while on the boundary of $C$ the trajectory should move on the sphere $x^2+y^2+(z+h)^2=1$ untill reaching the plane $z=0$ and then it moves on the intersection of the two spheres.
 
While in the interior of  $C$,  the control can change sign from $-1$ to $1$ or from $1$ to $-1$. Certainly, the control should be $1$ right before reaching the boundary and $-1$ right before arriving at $C_T$. Changes of the control  from $1$ to $-1$ or $-1$ to $1$ before reaching the boundary translate into time waste  and leads to smaller values of $x(T)$. 
It then follows that the optimal control should be of the form
\begin{equation}\label{opt:control}
u(t)=\left\{
\begin{array}{cl}
1,& t\in [0,\tilde{t}],\\[1mm]
-1, & t\in\ ]\tilde{t},T],
\end{array}
\right.
\end{equation}
 for some value   $\tilde{t}\in ]0,T[$.
\medskip

After the modification \eqref{tilde:psi},  the data of the problem satisfy the conditions under which Theorem \ref{main} holds. We now show that the conclusions of Theorem \ref{main} completly identify the structure \eqref{opt:control} of  the optimal control.

From  Theorem \ref{main} we deduce the existence of  $ \lambda \geq 0$, $ p,~q,~r \in BV([0,T],R)$, finite signed Randon measures 
 $ \eta_1$ and $\eta_2$,  null respectively in  $$I^{1}_b=\left\{ (x,y,z)\mid  x^2+y^2+(z+h)^2-1<0\right\}$$ and $$
  I^{2}_b=\left\{ (x,y,z)\mid x^2+y^2+(z-h)^2-1<0\right\},$$
$ \xi_{i}\in L^\infty([0,T],R)$, with $i=1,2  $, where $ \disp  \xi_{i}(t) \geq 0 \ \text{ a. e. } t$ and $\xi_{i}(t) = 0, \ t \in I^{i}_b,  $ 
such that 
$$\begin{array}{rl}
\text{(i)} & \begin{bmatrix}
\dot x(t)\\ \dot y(t)\\ \dot z(t)
\end{bmatrix}
=
\begin{bmatrix}
0 & \sigma & 0\\
0 & 0 & 0\\
0 & 0 & 0
\end{bmatrix}
\begin{bmatrix}
x\\ y\\ z
\end{bmatrix}
+\begin{bmatrix}
0\\ u\\ 0
\end{bmatrix}-2\xi_1\begin{bmatrix}
x\\ y\\ z+h
\end{bmatrix}-2\xi_2\begin{bmatrix}
x\\ y\\ z-h
\end{bmatrix}
\\[6mm]
\text{(ii)} & d\begin{bmatrix}
 p\\ q\\ r
\end{bmatrix} =  \begin{bmatrix}
0 & 0 & 0\\
-\sigma & 0 & 0\\
0 & 0 & 0
\end{bmatrix}\begin{bmatrix}
p\\ q\\ r
\end{bmatrix} dt
\\
& \hspace{2cm}
+2(\xi_1+\xi_2) \begin{bmatrix}
p\\ q\\ r
\end{bmatrix} dt
+2\begin{bmatrix}
x\\ y\\ z+h
\end{bmatrix} d\eta_1
+2\begin{bmatrix}
x\\ y\\ z-h
\end{bmatrix}d\eta_2,\\[6mm]
\text{(iii)} & \begin{bmatrix}
p\\ q\\ r
\end{bmatrix}(T) = 
\begin{bmatrix}
\lambda\\ 0\\ 0
\end{bmatrix}+\mu
\begin{bmatrix}
y_2\\ -\delta\\ 0
\end{bmatrix},\text{ where } \mu\geq 0,\\[5mm]
\text{(iv)} & \xi_1(xp+yq+(z+h)r)=0,\; \xi_2(xp+yq+(z-h)r)=0,\\[2mm]
\text{(v)}  &\text{the meaures } (xp+yq+(z+h)r)d\eta_1\text{ and }(xp+yq+(z-h)r)d\eta_2\\
& \hspace{0.5cm} \text{  are nonnegative,}
\\[2mm]
\text{(vi)} 
& \max_{u\in [-1,1]}uq=\hat{u}q.
\end{array}
$$
where  $\hat u$ is the optimal control.

Let $t_1$ be the instant of time when the trajectory reaches the  shere  $x^2+y^2+(z+h)^2=1$, 
 $t_2$ the instant of time when the trajectory reaches the intersection of the two spheres and $t_3$ be the instant of time the trajectory leaves the boundary of $C$. We have $0<t_1<t_2<t_3<T$.
 
 Next we show that the multiplier $q$ changes sign only once and so identifing the structure \eqref{opt:control} of the optimal control in a unique way.
 We start by looking at the case when    $t=T$. We have
$$\left[\begin{array}{c}
p\\ q
\end{array}\right] (T) = 
\left[\begin{array}{c}
\lambda\\ 0
\end{array}\right]+\mu
\left[\begin{array}{c}
y_2\\ -\delta\
\end{array}\right].$$
Starting from $t=T$, let us go  backwards in time  until the instant $t_3$  when the trajectory leaves the boundary of $C$.  If  $q(T)=0$, then $p(T)=\lambda>0$ and  we would have $q(t)>0$ for  $t \in ]t_3,T[$ (see (ii) above), which is impossible. We then  have $p(T)>0$ and $q(T)<0$ and, in $]t_3,T[$, since  $\sigma$ is small, the vector $(p(t),q(t))$ does not change much. At $t=t_3$, the vector $(p,q)$ has a jump and such jump can only occur along the  vector $(x(t_3),y(t_3))$.  Therefore, we have $p(t_3-0)>0$ and $q(t_3-0)<0$.

Let us now consider $t\in ]t_2,t_3[$. We have the following
\begin{enumerate}

\item when $ t\in [t_2,t_3]$, we have  $z=0$;

\item condition (i) above implies that  $\xi_1=\xi_2=\xi$, $\xi>0$ since, otherwise the motion along $x^2+y^2=1-h^2$ would not  be  possible;

\item from $0=\frac{d}{dt}(x^2+y^2)=\sigma 2xy-8\xi x^2+2uy-8\xi y^2$ we get  $\xi=\frac{\sigma xy+uy}{4(1-h^2)}$;

\item condition (iv) implies that  $r=0$  leading to  $xp+yq=0$. Since $x<0$, $y>0$, then $q=0$ implies $p=0$;

\item condition (ii) implies  $d\eta_1=d\eta_2=d\eta$;

\item $0=d(xp+yq)=uqdt+4(1-h^2)d\eta$ $\Rightarrow$ $\frac{d\eta}{dt}=-\frac{uq}{4(1-h^2)}$;

%\item $p>0,\; q<0$, $t>t_3$;
%
%\item (iv) $\Rightarrow$ $p>0$, $q>0$, $t=t_3$ (the jump is along $(x,y,0)$);

\item from the above analysis we deduce that
\begin{eqnarray*}
&& \dot{p}=
\frac{\sigma xy+uy}{(1-h^2)}~p-\frac{xuq}{(1-h^2)},\\
&& \dot{q}=-\sigma p+\frac{\sigma xy}{(1-h^2)}~q.
\end{eqnarray*}
Thus, $(p,q)$ is a solution to a linear system and it can never be equal to zero. It follows that $q$ cannot be zero because  $q=0$ implies $p=0$. Since   $q\neq 0$, we have  $q>0$.

\end{enumerate}

Let us consider the case when $t=t_2$. We  claim  that $$(p(t_2-0),q(t_2-0))\neq (0,0).$$ Seeking a contradiction, assume that it is  $(p(t_2-0),q(t_2-0))=(0,0)$. Then  we  have $$(p(t_2+0), q(t_2+0))=(0,0)+(2x_2(t_2), 2y_2(t_2))(d\eta_1+d\eta_2)$$ and such jump has to be   normal to  
$(x(t_2),y(t_2))$ since $r(t_2+0)=0$ (see (iv)).
It follows  that $(x^2(t_2)+y^2(t_2))(d\eta_1+ d\eta_2)=0$ and, since  $x^2(t_2)+y^2(t_2)>0$, we get   $d\eta_1+ d\eta_2=0$, proving our claim. 

\medskip

We now consider  $ t\in ]t_1,t_2[$. It is easy to see that $\xi_2=0$ and $d \eta_2=0$. We also deduce that 

\begin{enumerate}

\item $0=\frac{d}{dt}(x^2+y^2+(z+h)^2)=2\sigma xy+2uy-4\xi_1 y^2-4\xi_1 x^2-4\xi_1(z+h)^2$ which implies that $\xi_1=\frac{\sigma xy+uy}{2}$;

\item also $0=d(xp+yq+(z+h)r)=uqdt+2d\eta_1$ implies that  $\frac{d\eta_1}{dt}=-\frac{uq}{2}$;

\item from the above we deduce that 
\begin{eqnarray*}
&& \dot{p}=(\sigma xy+uy)p-xuq,\\
&& \dot{q}=-\sigma p+\sigma xyq.
\end{eqnarray*}
Thus $(p,q)$ is a solution to a linear system and never is equal to zero. Second equation implies that  if $q=0$ then $\dot{q}\neq 0$. Hence $q>0$.

\end{enumerate}

Now we need to consider  $t=t_1$. We   claim  that 
$$(p(t_1-0),q(t_1-0),r(t_1-0))\neq (0,0,0).$$ 
Let us then assume that it is  $(p(t_1-0),q(t_1-0),r(t_1-0))=(0,0,0)$.  It then follows that  $(p(t_1+0),q(t_1+0),r(t_1+0))=(0,0,0)+(2x(t_1)d\eta_1, 2y(t_1)d\eta_1, 2(z(t_1)+h)d\eta_1)$. We now show that there is no such jump. Set $r(t_1-0)=r_0$.  Then it follows from (iv) that 
 $(x(t_1)\cdot 0+y(t_1)\cdot 0 +(z(t_1)+h))r_0=0$ which  implies that $r_0=0$. 
 We also have   $(x^2(t_1)+y^2(t_1)+(z(t_1)+h)^2)d\eta_1=0$ from (v). But this implies that $d\eta_1=0$. Consequently,  the multipliers do not exhibit a jump at $t_1$. 

%It remains to study the behaviour of the multipliers while the  trajectory lives on the boundary of $C$. 
%There, the multiplier $q$ cannot be  $0$ for some interval of time. This is because if $q\equiv 0$, then 
%using $\rm (ii)$ we get
%$$
%0\leq \sigma pdt=y(d\eta_1+d\eta_2).
%$$
%But from $\rm (v)$ we have $d\eta_i\leq 0$, $i=1,2$. Hence $d\eta_i= 0$, $i=1,2$, and $p\equiv 0$ and $dp(\tilde{t})\leq 0$. This is impossible because $p(t)>0$, $t>\tilde{t}$. 

From the previous analysis we deduce that  $q$ should be positive almost everywhere on   the  boundary.   It then follows that 
to find the optimal solution we have to analyze admissible  trajectories with the controls with the structure \eqref{opt:control}
%$$
%u(t)=\left\{
%\begin{array}{cl}
%1,& t\in [0,\tilde{t}],\\
%-1, & t\in\ ]\tilde{t},T],
%\end{array}
%\right.
%$$
and choose the optimal value of $\tilde t$.

\begin{center}
\textbf{Acknowledgements}
\end{center}

The authors gratefully   thank the support of  Portuguese Foundation for Science and Technology
(FCT) in the framework of the Strategic Funding UIDB/04650\-/2020. 

Also we thank the support  by the ERDF - European Regional Development Fund through the Operational Programme for Competitiveness and Internationalisation - COMPETE 2020, INCO.2030, under the Portugal 2020 Partnership Agreement and by National Funds, Norte 2020, through CCDRN and FCT, within projects \textit{To Chair}  (POCI-01-0145-FEDER-028247), \textit{Upwind} (PTDC/EEI-AUT/31447/2017 - POCI-01-0145-FE\-DER\--03\-14\-47)
and \textit{Systec R\&D unit} (UIDB/00147/2020).

%References

\end{document}